\documentclass[a4paper,11pt]{amsart}

\usepackage{amsthm}
\usepackage[utf8]{inputenc}
\usepackage{verbatim}
\usepackage{amsmath}
\usepackage{amsfonts}
\usepackage{stackrel}
\usepackage{amssymb}
\usepackage{graphicx}
\usepackage[svgnames]{xcolor}
\usepackage{lipsum}
\usepackage{svg}
\usepackage{caption}
\usepackage{epigraph}
\usepackage{ebproof}
\usepackage[cal=boondoxo]{mathalfa}
\usepackage[all,cmtip]{xy}
\usepackage{mathtools}
\usepackage{color}
\usepackage{mathrsfs}
\usepackage{subfig}
\usepackage{framed}
\usepackage{hyperref}
\usepackage{tikz}
\usepackage{tikz-cd}

\usepackage[all,cmtip]{xy}

\usetikzlibrary{calc}
\usetikzlibrary{matrix,arrows}
\usepackage{tikz-cd}
\usepackage[square, sort, numbers]{natbib}

\usepackage{calligra} 

    {\endMakeFramed}

    {\endMakeFramed}

\makeatletter
\g@addto@macro\bfseries{\boldmath}
\makeatother

\theoremstyle{definition}

\theoremstyle{definition}

\newtheorem{rem}[equation]{Observation}

\theoremstyle{plain}
\newtheorem{thm}{Theorem}[section]
\newtheorem*{thm*}{Theorem}

\newtheorem{lem}[thm]{Lemma}

\newtheorem{defn}[equation]{Definition}

\newcommand\ck{\mathcal {K}}

\newcommand\cx{\mathcal {X}}

\title{Weak saturation and weak amalgamation property}
\author{Ivan Di Liberti}
\date{9 September 2017}

\address{
Ivan \textsc{Di Liberti}: \newline
Department of Mathematics and Statistics\newline
Masaryk University, Faculty of Sciences\newline
Kotl\'{a}\v{r}sk\'{a} 2, 611 37 Brno, Czech Republic\newline
\href{mailto:diliberti@math.muni.cz}{\sf diliberti@math.muni.cz}
}

\begin{document}
\maketitle

\begin{abstract}
The two model-theoretic concepts of \emph{weak saturation} and \emph{weak amalgamation} property are studied in the context of accessible categories. We relate these two concepts providing sufficient conditions for existence and uniqueness of weakly saturated objects of an accessible category $\ck$. We discuss the implications of this fact in classical model theory.
\end{abstract}

\tableofcontents

\section{Introduction}

In past years some work has been done in the direction of moving Fraïssé theory to a more category theoretic environment with two main purposes.  First, it looks like this environment is a natural habitat in which some argument become trivial and evident. Second, some results can be generalized from Fraïssé classes to well behaved categories. \newline

Two references for this work are Rosický \citep{Rsaturated} and Kubiś \citep{Kfrsequences} \citep{Kweak} in two very different manner. In the first one Rosický presents categorical aspects of saturation, moving existence and uniqueness of saturated objects to accessible categories. In the second Kubiś studies weakly $\omega$-saturated objects, again providing conditions to ensure existence and uniqueness. Here we bridge this two different approaches and hypotheses to obtain a fair generalization of both. \newline

A significant relaxing of the amalgamation property, called the \textit{weak amalgamation property}, was discovered by Ivanov \cite{IGenericEx} and independently by Kechris and Rosendal \cite{KSTurbulence} during their study of generic automorphisms in model theory.

\begin{defn}
A category $\ck$ has the weak amalgamation property if, for any object $A$ there is an arrow $A \stackrel{e}{\to} B$ such that any span like the following 

\begin{center}
\begin{tikzpicture}[scale=1.8]
\node (xa1) at (0,0) {$A$};
\node (xa2) at (0.5,0.5) {$B$};
\node (x1) at (0.5,1.5) {$C$};
\node (x2) at (1.5,0.5) {$D$};
\draw[->] (xa1) to node[below right]{$e$} (xa2);
\draw[->] (xa2) to node[left]{$f$} (x1);
\draw[->] (xa2) to node[below]{$g$} (x2);
\end{tikzpicture}
\end{center}

 can be completed to a square such that the diagram below is commutative.

\begin{center}
\begin{tikzpicture}[scale=1.8]
\node (xa2) at (0.5,0.5) {$A$};
\node (x1) at (0.5,1.5) {$C$};
\node (x2) at (1.5,0.5) {$D$};
\node (x3) at (1.5,1.5) {$E$};
\draw[->] (xa2) to node[right]{} (x1);
\draw[->] (xa2) to node[below]{} (x2);
\draw[->] (x2) to node[below]{} (x3);
\draw[->] (x1) to node[below]{} (x3);
\end{tikzpicture}
\end{center}

We call such an arrow $A \stackrel{e}{\to} B$ \textit{amalgamable}.

\end{defn}

This definition generalizes the \emph{cofinal amalgamation property}, which was already around in literature. This definition is given to present a framework for those classes where amalgamation property fails. Still one can pay a price to amalgamate over a model, embedding it in an \text{amalgamable hull}. \newline

Our aim is to convince the reader that the natural environment in which (weak) Fraïssé theory blooms is that of accessible categories with directed colimits. Here, some additional hypotheses guarantee the existence and the uniqueness of weakly saturated objects, a generalization of Fraïssé limits.  For any definition in the subject of accessible categories we direct the reader to \citep{book}. \newline

Let $T$ be a first order theory. The category Emb$(T)$, whose objects are models and morphisms are embeddings, is an accessible category with directed colimits; this is a perfect environment to look at model theory from the perspective of accessible categories. Eventually, a model $M$ in Emb$(T)$ has cardinality less than $\lambda$ if and only if it is $\lambda$-presentable. This is our dictionary when translating cardinality in the language of accessible categories. We will use $\lambda$-presentability as a replacement for \emph{cardinality is less than $\lambda$}. \newline

\begin{defn}
An object $K$ is weakly $\lambda$-saturated when for any arrow $A \to K$ where $A$ is $\lambda$-presentable there exists $A \to B$, with $B$ $\lambda$-presentable such that for any prolongation $A \to B \to C$ where $C$ is $\lambda$ -presentable there is an arrow $C \to K$ making the obvious diagram commutative.
\end{defn}

In the case of weak $\omega$-saturation Kubiś gave conditions for existence and uniqueness in \citep{Kweak}.
We will generalize Kubiś' result to weak $\lambda$-saturation using weaker hypothesis.

We state here  the main theorems we will prove in the sequel.

\begin{thm*}
Let $\ck$ be a $\lambda$-accessible category with the weak amalgamation property and directed colimits, satisfying a certain smallness condition,  then any object $K$ has a map $K \to M$ where $M$ is weakly $\lambda$-saturated.
\end{thm*}

\begin{thm*} Let $\ck$ be a $\lambda$-accessible category having  directed colimits and the joint embedding property. Then any two weakly $\lambda$-saturated, $\lambda^+$-presentable objects are isomorphic.
\end{thm*}

\begin{thm*} Let $\ck$ be a $\lambda$-accessible category having  directed colimits and the joint embedding property. A weakly $\lambda$-saturated, $\lambda^+$-presentable is weakly $\lambda$-homogeneous.
\end{thm*}

For readers more familiar with model-theoretic language, one can summarize the above results as follows, in the special context of categories of models of theories.

\begin{thm*}
Let $\ck$ be a class of models of a first order theory  with the weak amalgamation property and directed colimits, satifying a smallness condition and the joint embedding property, then

\begin{itemize}
\item any model embeds in a weakly $\lambda$-saturated one. 
\item any two weakly $\lambda$-saturated models of cardinality $\lambda$ are isomorphic.
\item weakly $\lambda$-saturated models of cardinality $\lambda$ are weakly $\lambda$-homogeneous.
\end{itemize} 
\end{thm*}

\section{Existence}

In this section we study accessible categories with the weak amalgamation property. Our aim is to find hypotheses to ensure the existence of weakly saturated objects. We will also discuss the implications of having such objects.

\subsection{A smallness condition}

\begin{defn}
We say that a category $\ck$ with the weak amalgamation property satisfies the smallness condition if, given a $\lambda$-presentable object $A$  and an amalgamable arrow $A \to M$, there exists a $\lambda$-presentable object $B$ and arrows $A \to B, B \to M$ such that $A \to B$ is amalgamable and the diagram below commutes.

\begin{center}
\begin{tikzpicture}[scale=1.8]
\node (xa2) at (0,0) {$A$};
\node (x1) at (0,-0.8) {$M$};
\node (x2) at (0.8,0) {$B$};

\draw[->] (xa2) to node[right]{} (x1);
\draw[->] (xa2) to node[below]{} (x2);
\draw[->] (x2) to node[below]{} (x1);
\end{tikzpicture}
\end{center}
\end{defn}

In short this condition ensures that it is not necessary to enlarge an object by much to obtain an amalgamation base. We call this property (SC).
Then in \citep{Kweak} any object is $\omega$-presentable, thus the condition is trivially verified.

\subsection{Existence}
\begin{thm} \label{ex}
Let $\ck$ be a $\lambda$-accessible category with the weak amalgamation property and directed colimits verifying (SC),  then any object $K$ has a map $K \to M$ where $M$ is weakly $\lambda$-saturated.

\begin{proof}
This is a small object argument, all in all. Let us call $\cx_K $  the set of all diagrams of the following shape 

\begin{center}
\begin{tikzpicture}[scale=1.8]
\node (xa2) at (0,0) {$B$};
\node (xa1) at (-1,0) {$A$};
\node (x1) at (1,0) {$C$};
\node (x2) at (0,1) {$K$};
\draw[->] (xa1) to node[above]{$e$} (xa2);
\draw[->] (xa2) to node[right]{$u$} (x2);
\draw[->] (xa2) to node[above]{$g$} (x1);
\end{tikzpicture}
\end{center}

where $e$ is amalgamable, $A,B,C$ are $\lambda$-presentables. We index this class with a cardinal $\eta = \text{card}(\cx_K)$

Now we build a chain $K_i$, $i < \eta$.

$K_0=K$.  \newline
$K_{i+1}$ is given by weak amalgamation property as follows

\begin{center}
\begin{tikzpicture}[scale=1.8]
\node (xa2) at (0,0) {$B_i$};
\node (xa1) at (-1,0) {$A_i$};
\node (x1) at (1,0) {$C_i$};
\node (x2) at (0,1) {$K_{i}$};
\node (x3) at (1,1) {$K_{i + 1}$};
\draw[->] (xa1) to node[above]{$e_i$} (xa2);
\draw[->] (xa2) to node[left]{$k_{0i \cdot u_i}$} (x2);
\draw[->] (xa2) to node[below]{$g_i$} (x1);
\draw[->] (x1) to node[below]{} (x3);
\draw[->] (x2) to node[above]{$k_{i,i+1}$} (x3);

\end{tikzpicture}
\end{center}

 For limit steps take directed colimit of the chain. \newline 
Call $K^*$ the directed colimit of of $K_i$, $i<\eta$.
In this way we built a morphism $K \to K^*$. If this is not an amalgamable arrow, we prolong it by an amalgamable arrow and take the resulting prolongation as our $K^*$. So we can assume that $K \to K^*$ is an amalgamable arrow.   \newline
Now call 
\begin{itemize}
\item $K^*= M_0$
\item  $M_{\alpha+1}= M_{\alpha}^*$,
\end{itemize}
We will call $m_{i,j}$ the map between $M_i$ and $M_j$.
After $\lambda$ steps, take the colimit and call it $M$. We call the colimit injection $M_i \stackrel{m_i}{\to} M$. We want to prove that $M$ is weakly saturated.

For a map $A \stackrel{s}{\to} M$ with $A$ $\lambda$-presentable, since $M$ is a $\lambda$-directed colimit and $A$ is $\lambda$-presentable, one can factor

\begin{center}
\begin{tikzpicture}[scale=1.8]
\node (xa2) at (0,0) {$A$};
\node (xa1) at (0,1) {$M$};
\node (x1) at (1,0) {$M_i$};
\draw[->] (xa2) to node[left]{$s$} (xa1);
\draw[->] (xa2) to node[above]{$u$} (x1);
\draw[->] (x1) to node[right]{$m_i$} (xa1);
\end{tikzpicture}
\end{center}

As we have the smallness condition, one can factor $ m_{i,i+1} \circ u $ as follows. This map is in fact amalgamable because $m_{i,j}$ is so.

\begin{center}
\begin{tikzpicture}[scale=1.8]
\node (xa2) at (0,0) {$A$};
\node (xa1) at (0,1) {$M$};
\node (x1) at (2,0) {$M_{i+1}$};
\node (x2) at (2,-1) {$B$};
\draw[->] (xa2) to node[left]{$s$} (xa1);
\draw[->] (xa2) to node[above]{$m_{i,i+1} \circ  u$} (x1);
\draw[->] (x1) to node[right]{$m_{i+1}$} (xa1);
\draw[->] (xa2) to node[below]{$e$} (x2);
\draw[->] (x2) to node[left]{} (x1);
\end{tikzpicture}
\end{center}

Now, we claim that for any prolongation $B \to C$ of $A \to B$ one can find an arrow $C \to M$ making the obvious triangle commutative. This is clear, because the following diagram will be in the set $\cx_{M_{i+1}}$:

\begin{center}
\begin{tikzpicture}[scale=1.8]
\node (xa2) at (0,0) {$B$};
\node (xa1) at (-1,0) {$A$};
\node (x1) at (1,0) {$C$};
\node (x2) at (0,1) {$M_{i+1}$};
\draw[->] (xa1) to node[above]{$e$} (xa2);
\draw[->] (xa2) to node[right]{$m_{i,i+1} \circ  u$} (x2);
\draw[->] (xa2) to node[above]{$g$} (x1);
\end{tikzpicture}
\end{center}

\end{proof}
\end{thm}

%
%
%

\subsection{A comment on the converse}

\subsubsection{Existence of weakly saturated objects implies...}
It would be natural to ask for a converse of Theorem \ref{ex}. In this direction much can be said.

\begin{rem} If any object $K$ has a map $K \to M$ where $M$ is weakly $\lambda$-saturated, then and object in the subcategory of $\lambda$-presentable objects has the weak amalgamation property.
\end{rem}

On the other hand we do not think that condition (SC) is implied by the existence of weakly saturated objects. It is possible, however, to get something similar to condition (SC).

\begin{lem} If any object $K$ has a map $K \to M$ where $M$ is weakly $\lambda$-saturated, then any amalgamable arrow $A \to K$, where $A$ is $\lambda$-presentable can fit in a diagram like the one below

\begin{center}
\begin{tikzpicture}[scale=1.8]
\node (xa2) at (0,0) {$A$};
\node (x1) at (0,-0.8) {$K$};
\node (x3) at (0.8,-0.8) {$N$};
\node (x2) at (0.8,0) {$B$};

\draw[->] (xa2) to node[right]{} (x1);
\draw[->] (xa2) to node[below]{} (x2);
\draw[->] (x2) to node[below]{} (x3);
\draw[->] (x1) to node[below]{} (x3);
\end{tikzpicture}
\end{center}

where the map $A \to B$ is amalgamable and $B$ is $\lambda$-presentable.

\end{lem}

(SC) is a stronger version of this, when one can choose $N=K$ and the bottom arrow to be the identity.
We have not yet managed to show that this weaker version implies the existence of weakly saturated objects.

\subsubsection{Looking for natural conditions that imply (SC)} It is quite hard to find natural conditions that imply (SC). A natural request could be that any object $K$ is a $\lambda$-pure subobject  $K \to M$ where $M$ is weakly $\lambda$-saturated. In fact:

\begin{lem}  If any object $K$ is a $\lambda$-pure subobject  $K \to M$ where $M$ is weakly $\lambda$-saturated then (SC) is verified.
\end{lem}

Although this is a terrible assumption, in fact it is easy to prove:

\begin{lem}  An object $K$ which is a $\lambda$-pure subobject  $K \to M$ of a weakly $\lambda$-saturated, is weakly $\lambda$-saturated. 
\end{lem}

This considerations could lead to the definition of weakly $\lambda$-pure morphism.

\begin{defn}
A morphism $K \to M$ is weakly $\lambda$-pure if for any diagram

\begin{center}
\begin{tikzpicture}[scale=1.8]
\node (xa2) at (0,0) {$A$};
\node (x1) at (0,-0.8) {$K$};
\node (x3) at (0.8,-0.8) {$M$};

\draw[->] (xa2) to node[right]{} (x1);
\draw[->] (x1) to node[below]{} (x3);

\end{tikzpicture}
\end{center}

there is an amalgamable arrow $A \to B$ such that, for any prolongation $B \to M$ making the square commutative, there is a morphism $B \to K$ making the upper triangle commutative.

\begin{center}
\begin{tikzpicture}[scale=1.8]
\node (xa2) at (0,0) {$A$};
\node (x1) at (0,-0.8) {$K$};
\node (x3) at (0.8,-0.8) {$M$};
\node (x2) at (0.8,0) {$B$};

\draw[->] (xa2) to node[right]{} (x1);
\draw[->] (xa2) to node[below]{} (x2);
\draw[->] (x2) to node[below]{} (x3);
\draw[->] (x1) to node[below]{} (x3);
\draw[dashed, ->] (x2) to node[below]{} (x1);
\end{tikzpicture}
\end{center}

\end{defn}

But this notion would turn any morphism into a weakly $\lambda$-pure one and thus does not seem interesting. \newline 

This ping pong between too strong and too weak conditions is quite discouraging to us in the direction of finding a condition which is necessary and sufficient for the existence of weakly saturated objects.

\section{Unicity}
In this section we discuss uniqueness of weakly saturated objects in accessible categories with the joint embedding property. In the end we introduce weak homogeneity and prove that weakly saturated objects are weakly homogeneous.

\begin{defn}
A category $\ck$ has the joint embedding property if any pair of objects $A, B$ have maps $A \to S$, $B \to S$ with same codomain.
\end{defn}

We will need a technical lemma. Its interest will be clear soon in the section. 

\begin{lem}\label{connectness}
If $\ck$ a $\lambda$-accessible category with the joint embedding property, than for any weakly $\lambda$-saturated object $K$ and any $\lambda$-presentable object $A$ there is an arrow $A \to K$.

\begin{proof}
Since the category is $\lambda$-accessible there is a $\lambda$-presentable  object $B$ and a map $B \to K$. Since $K$ is weakly saturated there is a map $B \to C$ whose prolongations can be extended. We now use joint embedding property on objects $A$ and $C$. This gives use two maps $C \to D$ and $A \to D$. $D$ has a map $D \to K$ by saturation and so by composition $A \to D \to K$ we get a map $A \to K$ as required.
\end{proof} 
\end{lem}

\begin{thm} \label{unique} Let $\ck$ be a $\lambda$-accessible category with directed colimits and the joint embedding property. Then any two weakly $\lambda$-saturated, $\lambda^+$-presentable objects are isomorphic.

\begin{proof}
Let $K$ and $L$ be $\lambda$-saturated and $\lambda^+$-presentable. It is well known that there are smooth chains $k_{i,j}: K_i \to K_j, i \leq j < \alpha$ and $l_{ij}: L_i \to L_j, i \leq j < \beta$ such that $\alpha, \beta \leq \lambda, K_i, i < \alpha, L_j, j < \beta$ are $\lambda$-presentable and $K = \text{colim}_{i < \alpha} K_{i}, L = \text{colim}_{i < \beta} L_{i}$. Let $k_i: K_i \to K, i < \alpha$ and $l_j: L_j \to L, j < \beta$ be colimit cocones. We may assume $\beta \leq \alpha$ We can also assume that $\alpha, \beta \in \{1, \lambda \}$. \newline

These sequences have a special property that we will state just for $K$. $\forall i < \alpha$ there exists $j \geq i$ such that forall $u: K_{j} \to D$ there is a map $g: D \to k$ such that $K_i = g \circ  u  \circ k_{i, j}$.  We may assume, by passing to subsequences if necessary, that $j = i +1$. In this case can the map $k_{i, i +1}$ is amalgamable.   \newline

We want to define cofinal subchains $\bar{k}_{ij}: \bar{K}_i \to \bar{K}_j, i \leq j < \alpha$ and $\bar{l}_{ij}: \bar{L}_i \to \bar{L}_j, i \leq j < \beta$ in $k_{ij}$ and $l_{ij}$, together with two family of maps $h_i, q_i$ such that 
$$h_i : \bar{K}_i \to \bar{L}_i $$
 $$q_i: \bar{L}_i \to \bar{K}_{i+1} $$
  $$q_i \circ h_i = \bar{k}_{i,i+1} $$
   $$h_{j+1}  \circ q_j = \bar{l}_{j, j+1}. $$ 

Thus colimit maps $h$ and $q$ will be mutually inverse between colimit objects. Since the subchains are cofinal, the colimit objects are precisely $K$ and $L$. \newline

[First Step]
Lemma \ref{connectness} proves that there is a map $K_1 \to L$, thus we can factor this map trough a $\lambda$-presentable, say $K_1 \stackrel{f_1}{\to} L_{j_0}\to L$  for a suitable $j_0$. By the special property of the sequence $K_i$, mentioned above, there is a map $L_{j_0+1} \stackrel{g_1}{ \to} K_{i_1} \to K$  such that $$g_1 \circ l_{j_0,j_0+1} \circ f_0 \circ k_{0,1} = k_{0, i_1}. $$

Set $$h_0 =   f_1 \circ k_{01}$$ $$ q_0 = g_1 \circ l_{i_0,i_1} $$ $$\bar{K}_0 = K_0$$ $$ \bar{L_0} = L_{j_0}. $$   \newline

[Isolated Step] 
We have  maps $h_i, q_i$ and we want to define $h_{i+1}, q_{i+1}$. Call $$\bar{K}_{i+1}=K_{i_t} = \text{cod}q_i. $$ By the special property of the sequence $L_j$ we have a map $ K_{i_t + 1} \stackrel{f_{i+1}  }{\to} L_{j_t} \to L$ and similarly we get $L_{j_t+1} \stackrel{g_{i+1}}{ \to} K_{i_m} \to K$.

Set $$h_{i+1} =   f_{i+1} \circ \bar{k}_{i_t,i_t+1}$$ $$ q_{i+1} = g_{i+1} \circ \bar{l}_{j_t,i_t +1} $$  $$ \bar{L}_{i+1} = L_{j_t}. $$

[Limit step] We take colimits. \newline

\end{proof}
\end{thm}

\subsection{Weak homogeneity}

As it is well known Fraïssé limits are $\omega$-homogenous. \citep{Kweak} obtains weak homogeneity for weakly $\omega$-saturated objects. Here we generalize this result to $\lambda$-saturated ones. In fact, the proof is the same of Kubiś.

\begin{defn} An object $M$ is weakly $\lambda$-homogeneous if for any amalgamable map $A \to B$ of  $\lambda$-presentable objects and any two prolongations as in diagram below

\begin{center}
\begin{tikzpicture}[scale=1.8]
\node (xa1) at (0,0) {$A$};
\node (xa2) at (0.7,0) {$B$};
\node (x1) at (1.4,0.7) {$M$};
\node (x2) at (1.4,-0.7) {$M$};
\draw[->] (xa1) to node[above]{$e$} (xa2);
\draw[->] (xa2) to node[above left]{$f$} (x1);
\draw[->] (xa2) to node[below left]{$g$} (x2);
\end{tikzpicture}
\end{center}

one can complete the triangle with an automorphism of $M$ in such a way that the diagram below is commutative 

\begin{center}
\begin{tikzpicture}[scale=1.8]
\node (xa2) at (0.7,0) {$A$};
\node (x1) at (1.4,0.7) {$M$};
\node (x2) at (1.4,-0.7) {$M$};
\draw[->] (xa2) to node[above left]{$f \circ e$} (x1);
\draw[->] (xa2) to node[below left]{$g \circ e$} (x2);
\draw[->] (x1) to node[below left]{} (x2);
\end{tikzpicture}
\end{center}

\end{defn}

\begin{thm}
Let $M$ be a weakly $\lambda$-saturated, $\lambda^+$-presentable object in an accessible category $\ck$ having  directed colimits and the joint embedding property, then $M$ is weakly $\lambda$-homogeneous.

\begin{proof}
Here there is not much to prove. One can follow the proof of the Theorem \ref{unique} and find the automorphism as a colimit map along the chain.
\end{proof}

\end{thm}

\newpage

\section*{Acknowledgement}
Author was supported by the GAČR project 17-27844S and RVO: 67985840.

\nocite{*}
\bibliography{bib.bib}
\bibliographystyle{alpha}

\end{document}